\documentclass[12pt,centertags,oneside, reqno]{amsart}
\usepackage{ifpdf}

\usepackage{amssymb}
\usepackage{braket}
\usepackage[all]{xy}
\usepackage{cancel}

\usepackage{amscd,amsxtra,calc}
\usepackage{cmmib57}
\usepackage{url}
\usepackage{ulem}
\usepackage{xcolor}
\usepackage[a4paper,width=16.2cm,top=3cm,bottom=3cm]{geometry}

\numberwithin{equation}{section}

\theoremstyle{plain}
\newtheorem{thm}{Theorem}[section]
\newtheorem{theorem}[thm]{Theorem}

\newtheorem{prop}[thm]{Proposition}
\theoremstyle{definition}

\newtheorem{dfn}[thm]{Definition}

\newtheorem{conjecture}[thm]{Conjecture}

\newtheorem{ques}[thm]{Question}

\numberwithin{equation}{section}

\title [Automorphisms of Calabi-Yau threefolds]{Automorphisms of Calabi-Yau threefolds from algebraic dynamics and the second Chern class}

\author{Keiji Oguiso}
\address{Mathematical Sciences, the University of Tokyo, Meguro Komaba 3-8-1, Tokyo, Japan, and National Center for Theoretical Sciences, Mathematics Division, National Taiwan University, 
Taipei, Taiwan}
\email{oguiso@ms.u-tokyo.ac.jp}

\thanks{The author is partially supported by JSPS Grant-in-Aid (A) 15H05738, JSPS Grant-in-Aid (B) 15H03611 and NCTS scholar program.}

\begin{document}

\maketitle

Throughout this report, $X$ is a normal projective variety of dimension $d$ over ${\mathbf C}$ and $f \in {\rm Bir}\, X$. This report consists of more speculations, rather than definite results, on automorphisms of strict Calabi-Yau threefolds in the view of the following question after Dinh, Sibony and D.-Q. Zhang (See \cite{Og14}, \cite{Og19} and \cite{Og23} and references therein for background, variants and relevant known results):
\begin{ques}\label{ques1} Find many examples of smooth $X$ admitting a {\rm primitive} biregular automorphisms of {\rm positive entropy}. (See below for definitions of some terms.)
\end{ques}
\begin{dfn}\label{dfn1} Let $H$ a very ample Cartier divisor on $X$. Then, the  $p$th dynamical degree of $f \in {\rm Bir}\, X$ is defined by
$$d_p(f) := \lim_{n \to \infty} ((f^{n})^{*}(H^{p}).H^{d-p})_{X}^{1/n} := \lim_{n \to \infty} (q_n^{*}(H^{p}).p_n^{*}(H^{d-p}))_{X_n}^{1/n} \ge 1,$$
where $p_n : X_n \to X$ is a birational morphism from a smooth projective variety $X_n$ such that $q_n : X_n \to X$ is a morphism and $f^n = q_n \circ p_n^{-1}$ (Hironaka's resolution). Let $\pi : X \dasharrow B$ be a dominant rational map to a normal projective variety $B$ with $f_B \in {\rm Bir}\, B$ such that $f_B \circ \pi = \pi \circ f$. In this case, the relative $p$th dynamical degree is defined by
$$d_p(f|\pi) = \lim_{n \to \infty} (f^{n}(H^{p}).H^{d-e}.\pi^{*}(H_B^e))_{X}^{1/n} 
\ge 1,$$
where $e := \dim\, B$ and $H_B$ is a very ample Cartier divisor on $B$.
\end{dfn}
Both $d_p(f)$ and $d_p(f|\pi)$ are well defined respectively by Dinh-Sibony \cite{DS05} and Dinh-Nguy$\hat{\rm e}$n \cite{DN11} (see also \cite{Tr20}). There are shown the existence of the limit, independence of $H$ and $H_B$, $d_{p+1}(f)d_{p-1}(f) \le d_{p}(f)^2$ (concavity), birational invariance of $d_p(f)$ and the product formula:
$d_p(f) := {\rm Max}_{k} \{d_k(f_B)d_{p-k}(f|\pi)\}$.
In particular $d_p(f) = d_p(f_B)$ if $e=d$, i.e., 
if $\pi$ is generically finite (thus recovers also the birational invariance). 

When $X$ is smooth and $f \in {\rm Aut}\, X$, the entropy of $f$ is ${\rm Max}_{p} \{\log d_p(f)\}$ (Gromov-Yomdin) and it is positive if and only if $d_1(f) > 1$ (by concavity), $(f^n)^{*} = (f^*)^n$ (obvious) and $d_p(f)$ is the spectral radius of $f^{*}|_{N^p(X)}$ also the one on $H^{p, p}(X, {\mathbf R})$ (an elementary fact on a linear map $f \in {\rm GL}(V, {\mathbf R})$ preserving a strictly convex closed cone $C \subset V$ with non-empty interior). 

\begin{dfn}\label{dfn2} $f \in {\rm Bir}\, X$ is imprimitive, if there are a dominant rational map $\pi : X \dasharrow B$ to a normal projective variety $B$ with $0 < e:= \dim\, B < d$ and $f_B \in {\rm Bir}\, B$ such that $\pi \circ f = f_B \circ \pi$. We call $f$ primitive if $f$ is not imprimitive.
\end{dfn}
This definition is due to De-Qi Zhang \cite{Zh09}. Using the product formula, we readily obtain the following useful criteria (See eg. \cite{OT15}, also \cite{KPS24} for other criteria):
\begin{prop} Let $f \in {\rm Bir}\, (X)$. Then $f$ is primitive if either $d_1(f) > 1$ and $\dim\, X =2$, $d_1(f) \not= d_2(f)$ and $\dim X = 3$, or $d_1(f) > d_2(f)$ (in any dimension).
\end{prop}
By \cite{Zh09} (See also \cite{Og14}), $X$ with primitive $f \in {\rm Bir}\, X$ is birational to either a smooth rationally connected variety, an abelian variety or a weak Calabi-Yau variety, i.e., a minimal projective variety with numerivally trivial canonical class and vanishing irregularity, under the assumption that minimal model program with abundance works in dimension $d$, which is expected to be true (and true if $d \le 3$).

When $d=2$, $X$ in Question \ref{ques1} is then birational to either ${\mathbf P}^2$, an abelian surface, a K3 surface or an Enriques surface. In each case, there found many interesting examples (See \cite{BC16}, \cite{OY20} and references therein). When $d=3$, a strict Calabi-Yau threefold (SCY3 for short), which is a {\it simply connected smooth projective threefold with trivial canonical line bundle}, is one of the most interesting tagert threefolds for Question \ref{ques1}. However, by now, we have only two such SCY3 (\cite{OT15} for $X_3$ and \cite[Lemma 4.6]{Og24} for $X_7$; it is easy to check $d_1(f) > d_2(f)$ for $f \in {\rm Aut}\, (X_7)$ there):
\begin{theorem} Two rigid SCY3 $X_k$ ($k =3$, $7$) have a primitive $f_k \in {\rm Aut}\, (X_k)$ of positive entropy. (See \cite{OS00} for explicit definitions of $X_k$).
\end{theorem}
We have the following intrinsic characterization of $X_k$ (\cite{OS00}):
\begin{dfn} Let $X$ be a SCY3. A surjective morphism $\varphi : X \to B$ to a normal projective variety $B$ with $\dim B > 0$ with connected fiber is called a $c_2$-contraction if $(c_2(X),\varphi^*H_B) = 0$ for an ample Caritier divisor $H_B$ on $B$, or equivalently, if $(c_2(X),\varphi^*N^1(B)) = 0$ (by the psuedo-effectivity of $c_2(X)$). Note that $c_2$-contraction $\varphi$ is never an isomorphism (as $c_2(X) \not= 0$ on $N^1(X)$ by the simply-connectedness of $X$). We call a $c_2$-contraction $\varphi_0 : X \to B_0$ maximal if $\varphi_0$ factors holomorphically through any $c_2$-contraction $\varphi : X \to B$. 
\end{dfn}
\begin{theorem}\label{thm1} If a SCY3 $X$ has a $c_2$-contraction, then $X$ has a maximal $c_2$-contraction $\varphi_0 : X \to B_0$. A SCY3 $X$ has a birational $c_2$-contraction $\varphi$ if and only if $X$ is isomorphic to one of $X_k$ ($k =3$, $7$). Also $\varphi$ is maximal in this case.
\end{theorem}
See \cite{OS00} for proof and details. Here partial vanishing of the second Chen class $c_2(X_k)$ on $\overline{{\rm Amp}\,(X_k)} \setminus \{0\}$ is natural for SCY3 $X$ of $|{\rm Aut}\, X| = \infty$:
\begin{prop} $c_2(X)^{\perp} \cap \overline{{\rm Amp}\,(X)}({\mathbf R}) \setminus \{0\} \not= \emptyset$ for a SCY3 $X$ of $|{\rm Aut}\, (X)| = \infty$.
\end{prop}
This is due to Wilson \cite{Wi97}. To my best knowledge, there is no known example such that $c_2(X)^{\perp} \cap \overline{{\rm Amp}\, (X)}({\mathbf R}) \setminus \{0\} \not= \emptyset$ but $c_2(X)^{\perp} \cap \overline{{\rm Amp}\,(X)}({\mathbf Q}) \setminus \{0\} = \emptyset$. So, it may be reasonable to ask:

\begin{ques}\label{ques2} Is $c_2(X)^{\perp} \cap \overline{{\rm Amp}\,(X)}({\mathbf Q}) \setminus \{0\} \not= \emptyset$ for a SCY3 $X$ of $|{\rm Aut}\, (X)| = \infty$? 
\end{ques}
There is also a notoriously difficult open conjecture (see eg. \cite{LOP18}):
\begin{conjecture}\label{conj1} Any nef integral divsior on a SCY3 is semi-ample.
\end{conjecture}
\begin{prop}\label{prop1} If both Question \ref{ques2} and Conjecture \ref{conj1} are affirmative, then there is no SCY3 $X$ with primitive 
$f \in {\rm Aut}\, X$ of positive entropy other than $X_k$ ($k= 3$ and $7$).
\end{prop}
Indeed, since $X$ has a primitive $f \in {\rm Aut}\, X$ of positive entropy, $X$ has a non-zero nef integral divisor $D$ with $(D.c_2(X)) = 0$ by Question \ref{ques2}. Then $|mD|$ for some $m >0$ gives a $c_2$-contraction $X \to B$ by Conjecture \ref{conj1}. Thus $X$ has a maximal $c_2$-contraction $\varphi_0 : X \to B_0$ by Theorem \ref{thm1}. By definition of maximality, we have a group homomorphism ${\rm Aut}\, X \to {\rm Aut}\, B_0$ being equivariant with respect to $\varphi_0$. Since $f \in {\rm Aut}\, X$ is primitive, it follows that $\dim\, B_0 =3$. Thus by Theorem \ref{thm1}, $X$ must be isomorphic to one of $X_k$ ($k=3$, $7$).
 
\medskip
\noindent

In the view of Proposition \ref{prop1}, it is very interesting to answer the following:
\begin{ques} Is there a SCY3 with primitive automorphism of positive entropy other than $X_k$ ($k= 3$ and $7$)?
\end{ques}

\end{document}